\documentclass[11pt,a4paper]{amsart}
\usepackage[english]{babel}
\usepackage{color}
\usepackage{amsmath,amssymb}
\usepackage{amsthm}
\usepackage[alphabetic, initials]{amsrefs}
\usepackage{geometry}
\usepackage{a4wide}

\usepackage{hyperref}
\usepackage{enumerate}

\renewcommand{\phi}{\varphi}
\renewcommand{\theta}{\vartheta}
\renewcommand{\epsilon}{\varepsilon}

\newcommand{\e}{\varepsilon}
\newcommand{\N}{\mathbb{N}}

\newcommand{\R}{\mathbb{R}}

\numberwithin{equation}{section}

\begin{document}

\title{Time-fractional Allen-Cahn equations \\
versus powers of the mean curvature}

\author[S. Dipierro]{Serena Dipierro}
\address{Serena Dipierro: Department of Mathematics and Statistics,
The University of Western Australia, 35 Stirling Highway, Crawley, Perth, WA 6009, Australia}
\email{serena.dipierro@uwa.edu.au}

\author[M. Novaga]{Matteo Novaga}
\address{Matteo Novaga: Dipartimento
di Mathematica, Universit\`a di Pisa, Largo Bruno Pontecorvo 5, 56127 Pisa, Italia}
\email{matteo.novaga@unipi.it.}

\author[E. Valdinoci]{Enrico Valdinoci}
\address{Enrico Valdinoci: Department of Mathematics and Statistics, The University of Western Australia, 35 Stirling Highway, Crawley, Perth, WA 6009, Australia}
\email{enrico.valdinoci@uwa.edu.au}


\thanks{EV was supported by the Australian Laureate Fellowship FL190100081 {\em Minimal surfaces, free
boundaries and partial differential equations}.
Part of this work has been carried out during a very pleasant visit of SD and EV to the University of Pisa, which we thank for the warm hospitality.}
\date{}

\dedicatory{}

\maketitle

\begin{abstract}
We show by a formal asymptotic expansion that level sets of solutions of a time-fractional Allen-Cahn equation evolve by
a geometric flow whose normal velocity is a positive power of the mean curvature.

This connection is quite intriguing, since 
the original equation is nonlocal and the evolution of its solutions depends on all previous states, but
the associated geometric flow is of purely local type, with no memory effect involved.
\end{abstract}

\section{Introduction}

Two very interesting, and apparently unrelated, topics have been intensively investigated in the contemporary mathematical literature, also in view of their applications and connections with other fields.\medskip

The first of these two topics focuses on the geometric flow of
hypersurfaces with a speed given by a positive power of the mean curvature
(when this power is equal to~$1$, the flow obviously reducing to the mean curvature flow~\cite{MR840401}).
The case of viscosity solutions has been treated in~\cite{MR1100211, MR1329522}.

The positive power of the mean curvature flow
is known to exist in case of closed and convex initial
data~\cites{MR2190140, MR2244700},
with finite-time extinction towards a point.

The problem has a very rich structure even in the plane,
where it provides surprising generalization of the curve shortening flow. For instance (see~\cites{MR845704, MR1888641, MR1949167})
when the power of the curvature is larger than~$\frac18$,
the only embedded homothetic solutions are circles, except when the power equals~$\frac13$, in which ellipses occur as well,
and when the power is below~$\frac18$ a new family of 
symmetric curves arises, resembling either circles or polygons.

The case of a volume-preserving flow has also been taken into account, see~\cite{MR4214340}.\medskip

The second topic of special interest is that of nonlocal equations,
and especially the time-fractional equations. This type of problems
has a classical flavor, dating back at least to the  tautochrone problem~\cite{Abel}, and emerges in several concrete examples,
such as viscoelastic fluids~\cite[Section~10.2]{MR1658022}
and diffusion along comb-like and fractal structures~\cites{Arkhincheev1991AnomalousDA, MR3856678, MR3967804} (see also
the introduction in~\cite{MR4321143} for several examples worked in full detail).

Among the several possible different choices of time-fractional derivatives, we recall the one introduced in~\cite{caputo}
in the context of dissipating models in geophysics and defined (up to a normalization constant that we omit for simplicity) as
$$ \partial^\alpha_t u(t):=\int_0^t \frac{\partial_t u(\tau)}{(t-\tau)^\alpha}\,d\tau,$$
with~$\alpha\in(0,1)$.

The technical advantage of this type of fractional derivative is
to often allow a consistent theory from initial conditions,
in the spirit of ordinary differential equations but comprising suitable ``memory effects''.
\medskip

To the best of our knowledge, these two topics,
namely, on the one side, geometric flows driven by powers of the mean curvature and, on the other side,
time-fractional equations driven by the so-called Caputo derivative
are considered as completely independent, and even quite separate
(though numerical simulations regarding classical and time-fractional Allen–Cahn models were performed in~\cite{MR4170335}, also leading to some conjectures about the dynamics of time-fractional models).

The goal of this paper is thus to show a deep link between
these two subjects, by considering the formal asymptotics
of a time-fractional Allen-Cahn
equation and by relating it with the geometric evolution
of level sets driven by powers of the mean curvature.

To this end, we will employ some asymptotic methods introduced
in~\cite{MR1245632} for the classical Allen-Cahn equation (see also~\cite{MR1067919} for formal asymptotics about the hyperbolic Allen-Cahn equation).
\medskip

More specifically, given~$\alpha\in(0,1)$ and~$\e>0$, we consider the time-fractional equation
\begin{equation}\label{pojqldwnfv-23erg}
\e^\alpha \partial^\alpha_t u=\e\Delta u+\frac{f(u)}{\e}. 
\end{equation}
In our setting, $u=u(x,t)$ with~$x\in\R^n$ and~$t\in[0,+\infty)$.

We will suppose that~$u$ takes values in~$[-1,1]$ and that~$f(u)$ is a bistable nonlinearity:
for concreteness, we focus on the case~$f(u):=u-u^3$.
We consider the global (strictly monotone) solution~$\gamma:\R\to[-1,1]$ of
\begin{equation}\label{Pga}
\begin{cases}
\gamma''+f(\gamma)=0 \quad {\mbox{in }}\R,\\
\gamma(0)=0,\\
\gamma(\pm\infty)=\pm1.\end{cases}
\end{equation}
With the above choice of~$f$, one has the explicit solution
\begin{equation}\label{TANHA}
\gamma(x):=\tanh\left(\frac{x}{\sqrt2}\right).\end{equation}

We consider the structural constants
\begin{equation}\label{PHI001} c_\alpha:=\int_{-\infty}^0\frac{\gamma'(\sigma)}{|\sigma|^\alpha}\,d\sigma\qquad{\mbox{and}}\qquad
C_\alpha:=\left(
\frac{(n-1)\gamma'(0)}{c_\alpha}\right)^{\frac1\alpha},
\end{equation}
and we take~$\phi_0$ satisfying
\begin{equation}\label{PHI0}
\begin{cases}
\dot\phi_0(t)=-\displaystyle\frac{C_\alpha}{\big(\phi_0(t)\big)^{\frac1\alpha}}\quad{\mbox{for }}t\in(0,+\infty),
\\ \phi_0(0)=1.\end{cases}
\end{equation}

{F}rom the geometric point of view, $\phi_0(t)$ describes the evolution in time of the radius of an $(n-1)$-dimensional sphere evolving
by the power~$\frac1\alpha$ of the mean curvature.

We consider the radial function
$$v_\e(r,t):=\gamma\left( \frac{r-\phi_0(t)}{\e}\right),$$
which, roughly speaking, models a layer function with spherical levels evolving by the power~$\frac1\alpha$ of the mean curvature
(as customary, here~$r=|x|$).
\medskip

The main result of this paper is that this radial layer function
is a solution (up to a small error) of the time-fractional Allen-Cahn equation, thus suggesting that level sets of solutions of the time-fractional Allen-Cahn equation have the tendency of evolving by powers of the mean curvature.

More precisely, via a formal expansion we will show that
\begin{equation}\label{TH1}
{\mbox{$v_\e$
solves~\eqref{pojqldwnfv-23erg} up to a small error in~$\e$.}}\end{equation}

\section{Asymptotic expansions}

We have that$$ \partial_tv_\e(r,t)=-\frac{\dot\phi_0(t)}{\e}\gamma'\left( \frac{r-\phi_0(t)}{\e}\right)$$
and, as a result,
$$ \partial^\alpha_t v_\e(x,t)=\int_0^t \frac{\partial_t v_\e(x,\tau)}{(t-\tau)^\alpha}\,d\tau=
-\int_0^t\frac{\dot\phi_0(\tau)}{\e\,(t-\tau)^\alpha}\gamma'\left( \frac{r-\phi_0(\tau)}{\e}\right)\,d\tau.
$$

Also, since, for every~$k\in\N$,
$$ \partial_r^k v_\e(r,t)=\frac1{\e^{k}}\partial^k\gamma\left( \frac{r-\phi_0(t)}{\e}\right),$$
we conclude that the Laplacian of~$v_\e$ can be written in the form
$$ \partial_r^2 v_\e(r,t)+\frac{n-1}r\partial_r v_\e(r,t)=
\frac1{\e^2}\gamma''\left( \frac{r-\phi_0(t)}{\e}\right)+\frac{n-1}{\e r}
\gamma'\left( \frac{r-\phi_0(t)}{\e}\right).$$

{F}rom these observations, we write the expression
$$
{\mathcal{E}}:=\e^\alpha
\partial^\alpha_t v_\e-\e\Delta v_\e-\frac{f(v_\e)}{\e}
$$ 
in the form
\begin{equation}\label{Gium}\begin{split}&
-\e^{\alpha-1}\int_0^t\frac{\dot\phi_0(\tau)}{(t-\tau)^\alpha}\gamma'\left( \frac{r-\phi_0(\tau)}{\e}\right)\,d\tau\\&\qquad-\frac1{\e}\gamma''\left( \frac{r-\phi_0(t)}{\e}\right)-\frac{n-1}{ r}
\gamma'\left( \frac{r-\phi_0(t)}{\e}\right)\\
&\qquad-\frac{1}\e f\left( \gamma\left( \frac{r-\phi_0(t)}{\e}\right)\right),
\end{split}\end{equation}
which, in view of the equation in~\eqref{Pga}, boils down to
\begin{equation}\label{Del00a}
-\e^{\alpha-1}\int_0^t\frac{\dot\phi_0(\tau)}{(t-\tau)^\alpha}\gamma'\left( \frac{r-\phi_0(\tau)}{\e}\right)\,d\tau-\frac{n-1}{ r}
\gamma'\left( \frac{r-\phi_0(t)}{\e}\right).
\end{equation}

We now denote by~${\mathcal{R}}_\e$ any quantity which is~$O(\e^k)$ for every~$k\in\N$: in particular, by~\eqref{TANHA},
\begin{equation}\label{Del00b} \gamma'\left( \frac{r-\phi_0(t)}{\e}\right)=\begin{cases} \gamma'(0) & {\mbox{ if }}r=\phi_0(t),\\
{\mathcal{R}}_\e& {\mbox{ if }}r\ne\phi_0(t).
\end{cases}\end{equation}

Moreover, by~\eqref{PHI0}, we know that~$\dot\phi_0<0$ and therefore
\begin{equation}\label{AKSxm-ORD}
{\mbox{$1=\phi_0(0)>\phi_0(\tau)>\phi_0(t)$ for all~$\tau\in(0,t)$.}}
\end{equation}

Now we recall the definition of~$c_\alpha$ in~\eqref{PHI001}, we pick~$\mu\in(0,1)$ and we claim that
\begin{equation}\label{Del00}
\int_0^t\frac{\dot\phi_0(\tau)}{(t-\tau)^\alpha}\gamma'\left( \frac{r-\phi_0(\tau)}{\e}\right)\,d\tau=\begin{cases}-c_\alpha\,\e^{1-\alpha}|\dot\phi_0(t)|^\alpha
+O(\e)&{\mbox{ if }}r=\phi_0(t),\\
O(\e^\mu) &{\mbox{ if }}r>\phi_0(t),\\
{\mathcal{R}}_\e& {\mbox{ if }}r<\phi_0(t).
\end{cases}
\end{equation}
Indeed, if~$r<\phi_0(t)$, we let~$a_0:=\phi_0(t)-r$ and we
have that~$r-\phi_0(\tau)<r-\phi_0(t)=-a_0$ for all~$\tau\in(0,t)$, due to~\eqref{AKSxm-ORD}, and accordingly, in light of~\eqref{TANHA},
$$ \gamma'\left( \frac{r-\phi_0(\tau)}{\e}\right)\le C\exp\left(-\frac{a_0}{\sqrt2\,\e}\right),$$
giving~\eqref{Del00} in this case.

If instead~$r>\phi_0(t)$, we argue as follows.
If~$|\phi_0(\tau)-r|\ge\e^\mu$, we have that~$\left|\frac{r-\phi_0(\tau)}\e\right|\ge\frac1{\e^{1-\mu}}$
and therefore
$$ \gamma'\left( \frac{r-\phi_0(\tau)}{\e}\right)={\mathcal{R}}_\e.$$
In this situation, we have that
\begin{equation}\label{mapP}
\int_0^t\frac{\dot\phi_0(\tau)}{(t-\tau)^\alpha}\gamma'\left( \frac{r-\phi_0(\tau)}{\e}\right)\,d\tau=
\int_{\mathcal{I}}\frac{\dot\phi_0(\tau)}{(t-\tau)^\alpha}\gamma'\left( \frac{r-\phi_0(\tau)}{\e}\right)\,d\tau+{\mathcal{R}}_\e,
\end{equation}
where
$$ {\mathcal{I}}:=\Big\{
\tau\in(0,t){\mbox{ s.t. }}|\phi_0(\tau)-r|<\e^\mu
\Big\}.$$
Now, the condition~$|\phi_0(\tau)-r|<\e^\mu$ boils down, as~$\e\searrow0$, to~$\tau=\phi_0^{-1}(r)<t$, and accordingly, for small~$\e$, we have that
$$ {\mathcal{I}}=\Big\{
\tau\in\R{\mbox{ s.t. }}|\phi_0(\tau)-r|<\e^\mu
\Big\}=\left( \phi_0^{-1}\big( r+\e^\mu\big),\phi_0^{-1}\big( r-\e^\mu\big)\right).$$
This and~\eqref{mapP} give that
\begin{eqnarray*}&&
\int_0^t\frac{\dot\phi_0(\tau)}{(t-\tau)^\alpha}\gamma'\left( \frac{r-\phi_0(\tau)}{\e}\right)\,d\tau=
\int_{\phi_0^{-1}( r+\e^\mu)}^{\phi_0^{-1}( r-\e^\mu)}\frac{\dot\phi_0(\tau)}{(t-\tau)^\alpha}\gamma'\left( \frac{r-\phi_0(\tau)}{\e}\right)\,d\tau+{\mathcal{R}}_\e\\
&&\qquad=\int_{\phi_0^{-1}( r+\e^\mu)}^{\phi_0^{-1}( r-\e^\mu)}\frac{\dot\phi_0\big(\
\phi_0^{-1}( r)+O(\e^\mu)
\big)}{\big(t-\phi_0^{-1}( r)+O(\e^\mu)\big)^\alpha}\gamma'\left( \frac{r-\phi_0(\tau)}{\e}\right)\,d\tau+{\mathcal{R}}_\e\\&&\qquad=\left(\frac{\dot\phi_0\big( \phi_0^{-1}( r)
\big)}{\big(t-\phi_0^{-1}( r)\big)^\alpha}+O(\e^\mu)\right)\int_{\phi_0^{-1}( r+\e^\mu)}^{\phi_0^{-1}( r-\e^\mu)}\gamma'\left( \frac{r-\phi_0(\tau)}{\e}\right)\,d\tau+{\mathcal{R}}_\e\\
&&\qquad=O(\e^\mu).
\end{eqnarray*}
This establishes~\eqref{Del00} in this case.

It remains to consider the case in which~$r=\phi_0(t)$. In this situation, given any~$\delta>0$,
to be taken conveniently small,
we know that, when~$\sigma\le-\frac\delta\e$,
$$ t-\phi_0^{-1}(\phi_0(t)-\e\sigma)\ge t-\phi_0^{-1}(\phi_0(t)+\delta)>0
$$
and we stress that the above quantity is bounded away from zero uniformly in~$\e$.

Therefore,
\begin{eqnarray*}&&
\int_0^t\frac{\dot\phi_0(\tau)}{(t-\tau)^\alpha}\gamma'\left( \frac{r-\phi_0(\tau)}{\e}\right)\,d\tau=
\int_0^t\frac{\dot\phi_0(\tau)}{(t-\tau)^\alpha}\gamma'\left( \frac{\phi_0(t)-\phi_0(\tau)}{\e}\right)\,d\tau\\&&\qquad=-\e\int_{\frac{\phi_0(t)-1}\e}^0\frac{\gamma'(\sigma)}{\big(t-\phi_0^{-1}(\phi_0(t)-\e\sigma)\big)^\alpha}\,d\sigma\\
\\&&\qquad=-\e\int_{-\frac{\delta}\e}^0\frac{\gamma'(\sigma)}{\big(t-\phi_0^{-1}(\phi_0(t)-\e\sigma)\big)^\alpha}\,d\sigma+O(\e).
\end{eqnarray*}
Now, when~$\e\sigma\in(-\delta,0)$, we let~$\psi_0:=\phi_0^{-1}$ and notice that
$$ \phi_0^{-1}(\phi_0(t)-\e\sigma)=
\psi_0(\phi_0(t)-\e\sigma)=t-\dot\psi_0(\phi_0(t))\e\sigma+O(\e^2\sigma^2)=
t-\frac{\e\sigma}{\dot\phi_0(t)}+O(\e^2\sigma^2)
$$
and, as a result,
\begin{eqnarray*}&&
\int_0^t\frac{\dot\phi_0(\tau)}{(t-\tau)^\alpha}\gamma'\left( \frac{r-\phi_0(\tau)}{\e}\right)\,d\tau=
-\e\int_{-\frac{\delta}\e}^0\frac{\gamma'(\sigma)}{\left(\frac{\e\sigma}{\dot\phi_0(t)}+O(\e^2\sigma^2)\right)^\alpha}\,d\sigma+O(\e)\\&&\qquad=-\e^{1-\alpha}\int_{-\frac{\delta}\e}^0\frac{\gamma'(\sigma)}{|\sigma|^\alpha\left(\frac{1}{|\dot\phi_0(t)|}+O(\e\sigma)\right)^\alpha}\,d\sigma+O(\e)\\&&\qquad=-\e^{1-\alpha}\int_{-\frac{\delta}\e}^0\frac{\gamma'(\sigma)}{|\sigma|^\alpha}
\Big(|\dot\phi_0(t)|^\alpha\big(1+O(\e\sigma)\big)\Big)
\,d\sigma+O(\e)\\&&\qquad=-\e^{1-\alpha}|\dot\phi_0(t)|^\alpha\int_{-\frac{\delta}\e}^0\frac{\gamma'(\sigma)}{|\sigma|^\alpha}
\,d\sigma+O(\e)\\&&\qquad=-c_\alpha\,\e^{1-\alpha}|\dot\phi_0(t)|^\alpha
+O(\e),
\end{eqnarray*}
which concludes the proof of~\eqref{Del00}.

Now, in the light of~\eqref{Del00a}, \eqref{Del00b} and \eqref{Del00},
we find that, when~$r\ne\phi_0(t)$,
$$ {\mathcal{E}}=O(\e^{\alpha-1+\mu}),$$
which is infinitesimal as long as we choose~$\mu\in(1-\alpha,1)$, and when~$r=\phi_0(t)$,
$$ {\mathcal{E}}=c_\alpha\,|\dot\phi_0(t)|^\alpha-\frac{(n-1)\gamma'(0)}{r}+O(\e^\alpha)=O(\e^\alpha),$$
owing to~\eqref{PHI0}.

This establishes~\eqref{TH1}, as desired.

\vfill

\end{document}